\documentclass[reqno]{amsart}
\usepackage{enumerate,amsmath,amssymb}
\usepackage{pstricks,pst-node,pst-coil,pst-plot} 
\usepackage{color}
\usepackage[implicit,naturalnames,hyperfootnotes=false,final]{hyperref}											
\usepackage{auto-eqlabel}
\usepackage{ifpdf}

\theoremstyle{plain}

\newtheorem{thm}{Theorem}[section]
\newtheorem{theorem}[thm]{Theorem}

\newtheorem{prop}[thm]{Proposition}
\newtheorem{proposition}[thm]{Proposition}

\theoremstyle{remark}

\newtheorem*{remark*}{Remark}
\newtheorem{remark}[thm]{Remark}

\theoremstyle{definition}

\newtheorem{definition}[thm]{Definition}

\newcounter{mnotecount}[section]

\newcommand{\definedas}{\mathrel{\raise.095ex\hbox{\rm :}\mkern-5.2mu=}}
\newcommand{\asdefined}{\mathrel{=\mkern-5.2mu}\raise.095ex\hbox{\rm :}\;}

\def\epsilon{{\varepsilon}}
\def\phi{{\varphi}}

\let\<\langle 
\let\>\rangle 

\DeclareMathOperator{\sh}{sh}
\DeclareMathOperator{\ch}{ch}
\DeclareMathOperator{\cone}{cone}

\newcommand{\bR}{\mathbb{R}}

\newcommand{\bS}{\mathbb{S}}

\newcommand{\bH}{\mathbb{H}}

\newcommand{\bU}{\mathbb{U}}
\newcommand{\bM}{\mathbb{M}}

\renewcommand{\hbar}{\overline{h}}

\DeclareMathOperator{\tr}{tr}
\DeclareMathOperator{\divg}{div}

\newcommand{\ric}{\mathrm{Ric}}

\newcommand{\scal}{\mathrm{Scal}}

\ifpdf\hypersetup{pdftitle= {On the center of mass of asymptotically hyperbolic initial data sets},									pdfauthor	= {Carla Cederbaum, Julien Cortier, and Anna Sakovich},
pdfsubject	= {Definition of the center of mass of asymptotically hyperbolic initial data sets with positive mass in asymptotically anti-de Sitter spacetimes},
pdfkeywords = {constant mean curvature foliation, AdS spacetime, hyperbolic space, center of mass, Hamiltonian charges, CMC-center of mass}}
\fi

\begin{document} 
\title{\sc{On the center of mass of asymptotically hyperbolic initial data sets}}
\author{Carla Cederbaum}
\address{Carla Cederbaum: Mathematics Department, Eberhard Karls Universit\"at T\"ubingen, Germany}
\email{cederbaum@math.uni-tuebingen.de}

\author{Julien Cortier}
\address{Julien Cortier: Max Planck Institute for Mathematics, Bonn, Germany}
\email{jcortier@mpim-bonn.mpg.de}

\author{Anna Sakovich}
\address{Anna Sakovich: Max Planck Institute for Gravitational Physics (Albert Einstein Institute), Potsdam-Golm, Germany}
\email{anna.sakovich@aei.mpg.de}

\thanks{The authors are grateful to the MSRI, Berkeley, to the FIM, ETH Z\"urich, to the AEI, Potsdam, and to the MPIM Bonn for hospitality and financial support during part of this work. Carla Cederbaum is indebted to the Baden-W\"urttemberg Stiftung for the financial support of this research project by the Eliteprogramme for Postdocs. This paper is typeset in \LaTeX\ with extra packages by Chr.~Nerz.}

\begin{abstract}
We define the (total) center of mass for suitably asymptotically hyperbolic time-slices of asymptotically anti-de Sitter spacetimes in general relativity. We do so in analogy to the picture that has been consolidated for the (total) center of mass of suitably asymptotically Euclidean time-slices of asymptotically Minkowskian spacetimes (isolated systems). In particular, we unite -- an altered version of -- the approach based on Hamiltonian charges with an approach based on CMC-foliations near infinity. The newly defined center of mass transforms appropriately under changes of the asymptotic coordinates and evolves in the direction of an appropriately defined linear momentum under the Einstein evolution equations.
\end{abstract}

\maketitle
\section{Introduction}

% !TEX root = CenterofMassAH.tex
The notion of the 'center of mass' of a physical system is one of the oldest and most fundamental concepts in mathematical physics and geometry. Understanding the position and motion of the center of mass of a system is often the first step towards gaining an understanding of the overall dynamics of the system. However, once one moves beyond classical mechanics, the concept of center of mass becomes increasingly complicated and needs to be re-defined. For example, in special relativity, the center of mass of a matter distribution intricately depends on the chosen observer, see M{\o}ller \cite{Moller}. This dependence will necessarily become more involved in general relativity, where we need to allude to an entire 'family' of observers.

Moreover, in order to account for black holes and other purely gravitational phenomena, it seems unreasonable to expect that the center of mass of, for example, an isolated gravitating system, will arise as a quantity solely defined by the matter components of the system. To the contrary, a reasonable definition of the (total) relativistic center of mass of a gravitational system will have to also apply to vacuum spacetimes such as the Schwarzschild or the Schwarzschild-anti-de Sitter (or Kottler) black holes. From these black hole examples, we also see that the center of mass cannot be described as an 'event' (point) in the spacetime manifold. Any notion of center of mass must thus necessarily be more abstract.

A very satisfying definition has been achieved for isolated systems described as suitably asymptotically Euclidean time-slices of asymptotically Minkow\-skian spacetimes, see below. In this paper, we will study suitably asymptotically hyperbolic time-slices of asymptotically anti-de Sitter spacetimes and give a definition of their (total) center of mass. We prove that the center of mass transforms appropriately under changes of the asymptotic coordinates and evolves in the direction of an appropriately defined linear momentum under the Einstein evolution equations. Our definition does \emph{not} coincide with the definition of (total) center of mass via Hamiltonian charges suggested by Chru\'sciel, Maerten, and Tod \cite{ChruscielMaertenTod}, but rather with their definition of energy-momentum vector or linear momentum. Our ideas rely on the picture that has been consolidated for the (total) center of mass for isolated systems:\\

\paragraph*{\emph{Isolated systems}} Three main approaches have been pursued to define the (total) center of mass of an \emph{isolated system}, that is, of a suitably asymptotically Euclidean time-slice of an asymptotically Minkowskian spacetime. First, the center of mass of an isolated system was defined via a Hamiltonian charges approach going back to Arnowitt, Deser, and Misner \cite{ADM}, Regge and Teitelboim \cite{RT}, Beig and \'O Murchadha \cite{BO}, and finally Szabados \cite{Szabados-old}, see \cite{Szabados} and Section~\ref{subsec:AEADM}. Subsequently, Huisken and Yau gave a geometric definition of the center of mass of a (time-symmetric\footnote{\emph{Time-symmetric} means that the second fundamental form of the embedding of the asymptotically Euclidean Riemannian time-slice into the asymptotically Minkowskian spacetime vanishes.}) isolated system \cite{HuiskenYau}, based on a foliation by stable constant mean curvature (CMC-)spheres in a neighborhood of infinity. This was motivated by an idea of Christodoulou and Yau \cite{ChristodoulouYau}. Huisken and Yau then 'coordinate' this 'abstract definition' of a center of mass by taking the limit of the Euclidean centers of the leaves of the CMC-foliation in the asymptotic end, see Section~\ref{subsec:AECMC} for more details and various generalizations.

Recently, Chen, Wang and Yau \cite{ChenWangYauQL,CWYEuclidean} gave a third definition of the (quasi-local and total) center of mass of an isolated system via optimal embeddings into Minkowski space. Under suitable decay conditions -- and after a $3+1$-decomposition of the Chen-Wang-Yau center --, all three (total) centers of mass are known to coincide, see Section~\ref{sect:AE}. The coordinatization of the center of mass of an isolated system can be interpreted as a point in $\bR^{3}$, the reference space which arises as the target space of the chosen asymptotically Euclidean chart near infinity. It transforms adequately under Euclidean coordinate transformations and it is known to evolve in the direction of linear (ADM-)momentum under the Einstein equations, see Section~\ref{subsec:AEprop} for more details.\\

\paragraph*{\emph{Hyperbolic systems}} In this paper, we study the (total) center of mass of suitably asymptotically hyperbolic time-slices of asymptotically anti-de Sitter (AdS) spacetimes. We will call such time-slices \emph{hyperbolic systems}, see Section~\ref{sect:AAdS} for the precise definition and the assumed asymptotic behavior. As in the case of isolated systems, one can take (at least) three different approaches to defining the (total) center of mass of a hyperbolic system. In analogy to the (abstract) CMC-center of mass defined by Huisken and Yau for isolated systems, Rigger \cite{Rigger}, and Neves and Tian \cite{NevesTian1,NevesTian2} showed that (time-symmetric\footnote{\emph{Time-symmetric} means that the second fundamental form of the embedding of the asymptotically hyperbolic Riemannian time-slice into the asymptotically anti de Sitter spacetime vanishes.}) hyperbolic systems possess a unique CMC-foliation near infinity which we will call the \emph{abstract} center of mass of the system. To do so, they assume that the asymptotic coordinate chart is 'balanced', see Section~\ref{Riemannian} for more details and references. Mazzeo and Pacard \cite{MazzeoPacard} also studied CMC-foliations near infinity of asymptotically hyperbolic Riemannian manifolds, however with non-spherical conformal infinity.

In Section~\ref{Riemannian}, we show that it is possible to coordinate this definition of center of mass in a way that is similar to the approach taken by Huisken and Yau. To do so, we extend the existence result for CMC-foliations of hyperbolic systems beyond balanced coordinates via hyperbolic boosts and define a new notion of 'coordinate center' of a hypersurface of any Riemannian manifold that can be 'nicely' isometrically embedded into Minkowski spacetime. We show that hyperbolic systems with balanced coordinates necessarily have their coordinate center of mass at the coordinate origin of hyperbolic space. From our definition, it is immediate that the CMC-center of mass transforms adequately under changes of the asymptotically hyperbolic chart near infinity, namely just like the center of the canonical slice of the Schwarzschild-anti-de Sitter (Kottler) spacetime.

On the other hand, it is also possible to define the (total) center of mass of a hyperbolic system as a Hamiltonian charge arising from Killing vector fields of the asymptotic background AdS spacetime, thus imitating the Hamiltonian charges definition of the center of mass for isolated systems described above. This has been studied by Chru\'sciel and Nagy \cite{ChruscielNagy} and physically interpreted by Chru\'sciel, Maerten, and Tod \cite{ChruscielMaertenTod}. However, in Section~\ref{Charges}, we find that our CMC-center of mass coincides with what they called the (Hamiltonian charge) energy-momentum vector (up to a normalizing factor) rather than with what they called the (Hamiltonian charge) center of mass. As in the case of isolated systems, the Hamiltonian charges approach thus provides an explicit asymptotic formula for computing the (CMC-)center of mass of a hyperbolic system, too. Using the Hamiltonian charges approach, we show that the center of mass evolves in direction of the linear momentum under the Einstein evolution equations with cosmological constant $\Lambda=-n(n-1)/2$, where $n$ is the dimension of the system, see Section~\ref{Evolution}. Here, the linear momentum of a hyperbolic system is defined as the '(Hamiltonian charge) center of mass' from \cite{ChruscielMaertenTod}.

The third approach -- using optimal embeddings into the AdS spacetime --, has been announced to work consistently by Chen, Wang, and Yau \cite{ChenWangYauhyperbolic}.
 
 \begin{remark*}
 We do \emph{not} expect our results to identically carry over to 'hyperboloidal systems', that is, to asymptotically hyperbolic time-slices of \emph{asymptotically Minkow\-skian} rather than \emph{asymptotically anti-de Sitter} spacetimes, as the role played by the second fundamental form conceptually differs in the two situations. 
 \end{remark*}

\subsection*{Notation and conventions} We will denote the Euclidean and hyperbolic spaces by $(\bR^{n},\delta)$ and $(\bH^{n},b)$, respectively. Both are special cases of $n$-dimensional Riemannian 'model spaces' which we will generally denote by $(\mathbb{M}^{n},h)$, see Section~\ref{Riemannian}. A general $n$-dimensional Riemannian manifold $(M^{n},g)$ will be called \emph{asymptotically Euclidean} or \emph{asymptotically hyperbolic} if there exists a smooth coordinate chart $\varphi:M^{n}\setminus K\to \bR^{n}\setminus \overline{B}$ or $\varphi:M^{n}\setminus K\to \bH^{n}\setminus \overline{B}$ such that $\varphi_{*}g-\delta$ or $\varphi_{*}g-b$ decays to $0$ suitably fast. Here, $K\subset M$ is a compact set and $B$ is a geodesic ball in the respective model space. As usual, if $(M^{n},g)$ arises as a time-slice of a spacetime $(L^{n+1},\mathfrak{g})$, we will call $(M^{n},g,k)$ an \emph{initial data set} or a \emph{system}. Here, $k$ denotes the second fundamental form. In contrast, whenever $k$ is insignificant or $k=0$, we will call $(M^{n},g)$ \emph{Riemannian} or \emph{time-symmetric}, respectively. All systems are assumed to have strictly positive mass, where the mass of a hyperbolic system is as defined in Section \ref{sect:balanced}.\\

All manifolds are assumed to be smooth. If not clear from context, we will indicate the metric $g$ with respect to which we are taking a divergence, a trace, a curvature tensor, or a covariant derivative by a lower right index $g$, for example $\divg_{g}$, $\tr_{g}$, $\ric_{g}$, $\nabla_{g}$, etc.

Closed hypersurfaces of $(M^{n},g)$ will be denoted by $\Sigma$, and we usually implicitly assume that $\Sigma\subset M^{n}\setminus K$ or in other words, $\Sigma$ lies in the domain of the asymptotic chart $\varphi$. The canonical metric on the unit sphere $\bS^{n-1}$ will be denoted by $\sigma$, the corresponding area measure will be denoted by $d\mu_{\sigma}$. Submanifolds $\Sigma\subset (M^{n},g)$ automatically inherit the induced metric from the ambient manifold; we will abuse notation and denote the induced area measure on $\Sigma$ by $d\mu_{g}$. This should help to prevent confusion when we are discussing several metrics on the ambient manifold.

Whenever needed, we will denote the Minkowski spacetime by $(\bR^{n,1},\eta)$. Both $(\bR^{n},\delta)$ and $(\bH^{n},b)$ can be isometrically embedded into $(\bR^{n,1},\eta)$, Euclidean space as a hyperplane with vanishing second fundamental form $k=0$ and hyperbolic space as an umbilic hyperboloid with second fundamental form $k=b$. We will denote these isometric embeddings by $I:(\mathbb{M}^{n},h)\hookrightarrow(\bR^{n,1},\eta)$. Moreover, $O(n,1)$ will denote the Lorentz group (the group of linear isometries of Minkowski space $(\bR^{n,1},\eta)$) and $SO_0(n,1)$ the \emph{restricted Lorentz group}, defined as the connected component of $O(n,1)$ which contains the identity; it is the group of direct isometries of $(\bH^{n},b)$.

All results concerning CMC-foliations will only apply to the physically relevant dimension $n=3$. Our definition of center of mass, the Hamiltonian charges approach, and the evolution result apply to all $n\geq3$.

We will use the $\mathcal{O}_{k}$-notation for tensors as an abbreviation for $C^{k}$-fall off. More specifically, let $r=\vert\vec{x}\vert$ denote the radial coordinate on $\bR^{n}$ and let $f:\bR^{+}\to\bR$ be a smooth function. Then for any smooth tensor field $T$ on $\bM^{n}$, we will write $T=\mathcal{O}_{k}(f(r))$ as $r\to\infty$ if there exists a constant $C>0$ independent of $r$ and a geodesic ball $B\subset\bM^{n}$ such that
\begin{align*}\labeleq{def:O-notation}
\vert\nabla^{\alpha}_{h} T\vert_{h} &\leq C \,\vert f^{\vert\alpha\vert}(r)\vert
\end{align*}
holds for all $x\in\bM^{n}\setminus \overline{B}$ and for all multi-indices $\alpha$ with $\vert\alpha\vert\leq k$. We will abuse this notation and also write $T=\mathcal{O}_{k}(f(r))$ as $r\to\infty$ when $T$ is a smoothly $r$-dependent family of tensors on $\bS^{n-1}$. In that case, $\nabla$ refers to both tangential derivatives along $(\bS^{n-1},\sigma)$ and partial $r$-derivatives. For isolated systems, $f$ will be an inverse power of $r$, while for hyperbolic systems, $f$ will be an exponential function $f(r)=e^{-lr}$ with $l\in\bR$ when we work in polar coordinates such that $b=dr^{2}+\sh^{2} r \,\sigma$. Here, $\sh$ denotes the hyperbolic sine function. Similarly, the hyperbolic cosine and arccosine will be denoted by $\ch$ and $ \mathrm{arcch}$, respectively.

Finally, constant points/vectors in $\bR^{n}$ such as the coordinate center of mass $\vec{z}$ of an isolated system will be denoted by lower case letters with an arrow. Similarly, constant points in $\bH^{n}$ such as the coordinate center of mass of a hyperbolic system $\mathbf{z}$ will be denoted by bold lower case letters. Constant vectors in $\bR^{n,1}$ such as the 'mass vector' $\mathbf{P}$ will be denoted by bold upper case letters.

\subsection*{Structure of the paper} This paper is structured as follows: In Section~\ref{sect:AE}, we give a concise overview of the definition of the (total) center of mass for asymptotically Euclidean systems and related results. In Section~\ref{sect:balanced}, we summarize and extend known results on asymptotic CMC-foliations of asymptotically hyperbolic systems and give a definition of an associated hyperbolic coordinate center of mass. More general results on this CMC-center of mass are described and proven in Section~\ref{sect:CoM}. In Section~\ref{Charges}, we generalize our definition to asymptotically hyperbolic initial data sets and link it to the Hamiltonian charges approach. Finally, in Section~\ref{Evolution}, we prove that the center of mass evolves appropriately under the Einstein evolution equations.

\subsection*{Acknowledgements} We wish to thank Michael Eichmair and Romain Gicquaud for helpful and stimulating discussions.
\section{Summary of results in the asymptotically Minkowskian setting}\label{sect:AE}

% !TEX root = CenterofMassAH.tex
Three main approaches have been pursued to define the (total) center of mass of a suitably asymptotically Euclidean time-slice of an asymptotically Minkowskian spacetime. First, the center of mass of an isolated system was defined via a Hamiltonian charges approach, see Section~\ref{subsec:AEADM}. Subsequently, Huisken and Yau \cite{HuiskenYau} gave a geometric definition of the center of mass in the Riemannian context, based on a foliation by stable constant mean curvature (CMC-)spheres in a neighborhood of infinity, see Section~\ref{subsec:AECMC}. Both definitions are known to coincide, see Section~\ref{subsec:AEprop}.

Recently, Chen, Wang and Yau \cite{ChenWangYauQL,CWYEuclidean} gave a third definition of the (quasi-local and total) center of mass of an isolated system via optimal embeddings into Minkowski space. Under suitable decay conditions, their total center of mass definition coincides with the definition via Hamiltonian charges and CMC-foliations, see~\cite{CWYEuclidean}. We will not further discuss their center of mass definition as we do not intend to generalize it to the asymptotically anti-de Sitter setting. However, this generalization has been announced to work consistently by Chen, Wang, and Yau~\cite{ChenWangYauhyperbolic}.

We will not attempt to cite the weakest possible decay conditions and the strongest results, here. We rather focus on recalling the general ideas from the asymptotically Euclidean context that we will rely on in the asymptotically hyperbolic context in the following sections. Precise statements and a more complete list of references can be found in \cite{Cederbaum-Nerz}. 

First, let us recall the precise definition of asymptotically Euclidean and asymptotically Schwarz\-schil\-dean Riemannian manifolds and of asymptotically Euclidean initial data sets (isolated systems):

\begin{definition}
A Riemannian $n$-manifold $(M,g)$ is \emph{$(l,\tau)$-asymp\-to\-ti\-cally Euclidean (with one end)} if there exists a compact set $K \subset M$, a constant $R > 0$, and a diffeomorphism $\phi : M\setminus K \rightarrow \bR^{n}\setminus \overline{B}_R$ such that the metric $g^\phi \definedas \phi_* g$ satisfies
\begin{align*}\labeleq{eqEucl}
g^\phi &= \delta+ \mathcal{O}_{l}\left(\frac{1}{r^{\tau}}\right)
\end{align*}
as $r\to\infty$ and \emph{$(l,\tau)$-asymptotically Schwarzschildean (with one end)} if 
\begin{align*}\labeleq{eq:Schwarz}
g^\phi &= \left(1+\frac{m}{2r}\right)^{4}\delta+ \mathcal{O}_{l}\left(\frac{1}{r^{1+\tau}}\right)
\end{align*}
as $r\to\infty$. Here, $m\in\bR$ and we assume $\tau>0$ and $l\geq0$. Similarly if $(M,g,k)$ is an \emph{$(l,\tau)$-asymptotically Euclidean initial data set (with one end)} if $(M,g)$ is $(l,\tau)$-asymptotically Euclidean with respect to a diffeomorphism $\varphi$ and
\begin{align*}
\varphi_{*}k &=\mathcal{O}_{l-1}\left(\frac{1}{r^{2+\tau}}\right)
\end{align*}
as $r\to\infty$.
\end{definition}

\begin{remark}
Alternatively, instead of using $\mathcal{O}$-notation (or weighted $C^{l}$ spaces), one can formulate the asymptotic decay in terms of weighted Sobolev spaces. 
\end{remark}

The ADM-mass $m$ of a $3$-dimensional asymptotically Euclidean Riemannian manifold $(M^{3},g)$ is defined as follows \cite{ADM}:
\begin{definition}[ADM-mass]
The \emph{ADM-mass} $m$ of an $(l,\tau)$-asymp\-to\-ti\-cally Euclidean Riemannian $3$-manifold $(M^{3},g)$ with diffeomorphism $\varphi$ is defined as
\begin{align*}
 m \definedas{}& \frac{1}{16\pi}\lim_{r\to\infty}\sum_{i,j=1}^3\ \int_{\bS^2_r}(\partial_i{g^{\varphi}_{jj}} - \partial_i{g^{\varphi}_{ij}})\,\frac{x_i}r \,d\mu_{\delta}.\labeleq{mADM} 
\end{align*}
\end{definition}
This expression is well-defined and independent of the chart $\varphi$ if $\tau>1/2$, $l\geq2$, and if the scalar curvature $R$ is $L^{1}$-integrable, see Bartnik \cite{Bartnik} and Chru\'sciel~\cite{Chrus2}.

\subsection{The center of mass of an isolated system via Hamiltonian charges}\label{subsec:AEADM}
Using methods from the theory of Hamiltonian systems, Regge and Teitelboim \cite{RT} and Beig and {\`O} Murchadha \cite{BO} have constructed a Hamiltonian charge that can be interpreted as the center of mass of suitably asymptotically Euclidean manifolds with positive ADM-mass $m$. This has been generalized to suitably asymptotically Euclidean initial data sets by Szabados \cite{Szabados-old,Szabados}. Going beyond the asymptotically Euclidean setting, the asymptotic invariants approach to defining the center of mass (and the mass) of an $n$-dimensional Riemannian manifold with positive mass has been carried out by Michel \cite{MichelMass}, from where we take the following definition and theorem:

\begin{theorem}[Michel]
 Let $n \geq 3$ and let $(M,g)$ be an $n$-dimensional $(l,\tau)$-asymptotically Euclidean manifold of order $\tau > \frac{n-1}{2}$ and $l\geq2$ with respect to a diffeomorphism $\varphi$. Assume that the ($n$-dimensional) ADM-mass $m$ of $g$ is positive. Then the following limits exist and are finite:
 \begin{align*}\labeleq{def:hamiltonianCoM}
  z^i\definedas  \frac{1}{2m(n-1)\, \omega _{n-1}}\lim_{r \to \infty} \int _{\bS^{n-1}_{r}} \bU(x^i,e)(\nu _\delta)\,d\mu_{\delta},
 \end{align*}
where $e \definedas g - \delta$, and $\nu _\delta$ is the outward unit normal to the coordinate sphere $\bS^{n-1}_{r}$. The $1$-form $\bU(x^i,e)$ appearing in the integrand is given by the expression
\begin{align*}
 \bU(x^i,e) \definedas x^i \left(\divg_{\delta} e - d \tr_{\delta} e\right) - \iota _{\nabla_{\delta} x^i} e + \tr_{\delta} e\ dx^i.
\end{align*}
The vector $\vec{z}$ with components $z^i$ is called the \emph{(Hamiltonian) center of mass} of $(M,g)$ in the coordinate system induced by $\varphi$.
\end{theorem}

\begin{remark}
It is well-known that the expression \eqref{def:hamiltonianCoM} also converges if $\tau>\frac{n-2}{2}$ and $(M,g)$ is asymptotically even with respect to $\varphi$ (or, in other words, satisfies the Regge-Teitelboim conditions).
\end{remark}

\begin{remark}
By construction, the Hamiltonian center of mass transforms adequately under Euclidean motions on $\bR^{n}$ (in the image of $\varphi$). It has been shown under suitable decay conditions on the full initial data set\footnote{See Section \ref{Charges} for the definition of initial data sets.} $(M,g,k,\mu,J)$ and the lapse and shift of the evolution that the Hamiltonian center evolves under the Einstein equations such that $\frac{d}{dt}\left(m\vec{z}\right)=\vec{P}$, where $\vec{P}$ is the ADM-linear momentum, see~\cite{Szabados} and the references cited therein.
\end{remark}

\subsection{The center of mass of an isolated system via CMC-foliations}\label{subsec:AECMC}
Several authors define the center of mass of a suitably asymptotically Euclidean Riemannian manifold $(M,g)$ with ADM-mass $m>0$ as a foliation near infinity. Following Cederbaum and Nerz \cite{Cederbaum-Nerz}, we will call such definitions \emph{abstract} to contrast what we will call \emph{coordinate definitions of center of mass}, see below. The first and most important such abstract definition was given by Huisken and Yau \cite{HuiskenYau}, who defined the (abstract) \emph{CMC-center of mass} to be the unique foliation near infinity by closed, stable surfaces with constant mean curvature. This was motivated by an idea of Christodoulou and Yau \cite{ChristodoulouYau}. Huisken and Yau \cite[Theorem 4.2]{HuiskenYau} show that the CMC-foliation exists and is unique whenever the Riemannian manifold $(M^{3},g)$ is $(4,1)$-asymptotically Schwarzschildean and has positive mass. Their existence and uniqueness result has been generalized significantly; for an overview over other definitions of abstract center of mass notions and for a discussion of generalizations and optimal decay, please see \cite{Cederbaum-Nerz}.

In order to assign a center of mass vector $\vec{z}$ to an abstract center of mass given by a foliation near infinity, one can utilize the idea of a Euclidean coordinate center: In Euclidean geometry, any closed surface $\Sigma\hookrightarrow\bR^n$ has a \emph{Euclidean coordinate center} $\vec{z}_{\Sigma}$ defined by
\begin{align*}
\vec{z}_{\Sigma}\definedas \frac{1}{\vert\Sigma\vert} \int_{\Sigma} \vec{x} \,d\sigma.\labeleq{eucl_center_of_surface} \end{align*}
Now pick a fixed diffeomorphism $\varphi$ near infinity with respect to which a given Riemannian manifold $(M^{n},g)$ is asymptotically Euclidean. Then $\varphi$ allows us to define the \emph{Euclidean coordinate center $\vec{z}_{\Sigma;\varphi}$} of a closed surface $\Sigma\hookrightarrow M\setminus K$ via
\begin{align*}
\vec{z}_{\Sigma;\varphi}&\definedas \frac{1}{\vert\varphi(\Sigma)\vert_{\delta}}\int_{\varphi(\Sigma)}\vec{x}\,d\mu_{\delta}.
\end{align*}

Huisken and Yau \cite[Theorem 4.2]{HuiskenYau} subsequently define the Euclidean center of the CMC-foliation near infinity as the limit of the Euclidean coordinate centers $\vec{z}_{\Sigma;\varphi}$ outward along the foliation. More abstractly, this can be described by making the following definition, cited from \cite{Cederbaum-Nerz}.

\begin{definition}[Coordinate center of a foliation]\label{folicoordCoM}
Let $\{\Sigma_{\sigma}\}_{\sigma\ge\sigma_0}$ be a foliation near infinity of an $(l,\tau)$-asymptotically Euclidean Riemannian $n$-manifold $(M,g)$ with asymptotic coordinate chart $\varphi$ where $\sigma\geq\sigma_{0}$ is such that $\sigma\to\infty$ corresponds to enumerating the surfaces $\Sigma_{\sigma}$ outward to infinity. Let $\vec{z}_{\sigma;\varphi}$ denote the Euclidean coordinate center of the leaf $\Sigma_{\sigma}$. In case the limit exists, the \emph{coordinate center $\vec{z}_{\varphi}$ of the foliation $\{\Sigma_\sigma\}_{\sigma\geq\sigma_0}$} is given by 
\begin{align*}
\vec{z}_{\varphi}\definedas\lim_{\sigma\to\infty}\vec{z}_{\sigma;\varphi}.
\end{align*}
\end{definition}

When $\{\Sigma_{\sigma}\}_{\sigma\ge\sigma_0}$ is the CMC-foliation constructed by Huisken and Yau \cite[Theorem 4.2]{HuiskenYau}, we will call the induced coordinate center $\vec{z}_{\varphi}$ the (coordinate) CMC-center of mass of the slice and denote it by $\vec{z}_{\text{CMC}}$. Cederbaum and Nerz \cite{Cederbaum-Nerz} gave explicit examples demonstrating that the coordinate CMC-center of mass $\vec{z}_{\text{CMC}}$ need not converge in a $(4,1)$-asymptotically Schwarzschildean Riemannian $3$-manifold as studied by Huisken and Yau because the limit in Definition \ref{folicoordCoM} does not necessarily converge. However, it can be demonstrated that the limit and thus the coordinate CMC-center of mass $\vec{z}_{\text{CMC}}$ does converge if the manifold is $(4,1+\varepsilon)$-asymptotically Schwarzschildean Riemannian for any $\varepsilon>0$.

\subsection{Properties of the (total) center of mass definitions}\label{subsec:AEprop}
It is well-known that the coordinate CMC-center of mass of a given, suitably asymptotically Euclidean Riemannian $3$-manifold $(M^{3},g)$ coincides with its ADM-center of mass. This has been established by Huang \cite{Lan_Hsuan_Huang__Foliations_by_Stable_Spheres_with_Constant_Mean_Curvature}, Eichmair and Metzger \cite{EichmairMetzgerInventiones}, and Nerz \cite[Corollary~3.8]{nerz2013timeevolutionofCMC} under different assumptions on the asymptotic decay. Moreover, it has been demonstrated by Cederbaum that the coordinate CMC- and ADM-centers of mass converge to the Newtonian one in the Newtonian limit $c\to\infty$ \cite[Chap.~4,\;6]{Cederbaum_newtonian_limit} whenever the manifold is (vacuum) static. This further justifies the ADM- and CMC-definitions of center of mass. 

It is straightforward to demonstrate that both the coordinate CMC-center of mass $\vec{z}_{\text{CMC}}$ and the ADM-center of mass $\vec{z}_{\text{ADM}}$ transform adequately under asymptotic Euclidean motions performed in the image of the coordinate chart $\varphi$. The CMC-foliation itself however is independent of the chart $\varphi$ and thus does not get affected by changes of the coordinates. For more details on its existence and uniqueness without reference to a given asymptotic chart $\varphi$, see Nerz \cite{Nerz-characterization}.

\section{The CMC-center of mass of asymptotically hyperbolic manifolds}\label{sect:AAdS}

\subsection{In balanced coordinates}\label{sect:balanced}
% !TEX root = CenterofMassAH.tex
Before giving a summary of the results on the existence and uniqueness of CMC-foliations of asymptotically hyperbolic (Riemannian) manifolds, let us first recall some definitions and specify the exact fall-off conditions we will work with. In the sequel, we will choose a point in $\bH^n$ as the origin. Then in polar coordinates around this point, the hyperbolic metric is given by \begin{align*}\labeleq{eqHyperbolic}
b= dr^2 + \sh^2 r \,\sigma,
\end{align*}
where $r$ is the Euclidean distance to the origin and $r^2 = (x^1)^2 + \cdots +(x^n)^2$. We also agree to denote the inner product and the norm induced by $b$ on natural tensor bundles over $\bH^n$ by $\langle \cdot,\cdot \rangle$ and $\vert\cdot\vert$, respectively. The following definition of an asymptotically hyperbolic (Riemannian) manifold is a slight generalization of X.~Wang \cite[Definition 2.3]{WangMass}. 

\begin{definition}\label{defWang}
 A Riemannian $n$-manifold $(M,g)$ is \emph{$l$-asymptotically hyperbolic (with one end)} if there exists a compact set $K \subset M$, a constant $R > 0$ and a diffeomorphism $\phi : M\setminus K \rightarrow \bH^{n}\setminus \overline{B}_R$ such that the metric $g^\phi \definedas \phi_* g$ can be written in polar coordinates as
\begin{align*}\labeleq{eqWang}
g^\phi = dr^2 + \sh^2 r\, h_r ,
\end{align*}
where $h_r$ is a family of symmetric 2-tensor fields on $\bS ^{n-1}$ admitting the expansion
\begin{align*}\labeleq{eq:hr}
 h_r = \sigma + \mathbf{m}\, e^{-nr} + \mathcal{O}_{l}(e^{-(n+1)r})
\end{align*}
as $r\to\infty$. Here, $\mathbf{m}$ is a symmetric 2-tensor field on $\bS^{n-1}$, and we assume $l \geq2$.
\end{definition}

In particular, we recover the hyperbolic metric from \eqref{eqWang} when $h_r = \sigma$ for all $r>0$. The symmetric 2-tensor $\mathbf{m}$ is called the \emph{mass aspect tensor}, its trace $\tr _{\sigma} \mathbf{m}$ is the so-called \emph{mass aspect function}. Now let $\mathring{x}^i$,  $i=1,\ldots, n$, denote the restrictions of the Euclidean coordinate functions $x^i$ to $\bS^{n-1} \subset \bR^n$. The vector
\begin{align*}\labeleq{eqMassVectorRiem}
 \mathbf{P} \definedas \left(\int_{\bS^{n-1}} \tr _{\sigma} \mathbf{m}\,d\mu_\sigma, 
 \int_{\bS^{n-1}} \mathring{x}^1 \tr _{\sigma} \mathbf{m}\, d\mu_\sigma, \ldots, 
 \int_{\bS^{n-1}} \mathring{x}^n \tr _{\sigma} \mathbf{m}\, d\mu_\sigma\right)\in\bR^{n,1}
\end{align*}
will be referred to as the \emph{mass vector}, although the term \emph{energy-momentum vector} -- introduced in analogy to the asymptotically Euclidean setting -- is quite common in the literature, see for example X.~Wang \cite{WangMass}. The mass aspect tensor $\mathbf{m}$ and thus the components of the {mass vector} $\mathbf{P}$ clearly depend on the choice of the coordinate chart $\varphi$. 
However, $-\vert\mathbf{P}\vert_{\eta}^{2}$ is an \emph{asymptotic invariant}, meaning that it does not depend on the choice of the diffeomorphism $\varphi$ in Definition \ref{defWang}. The proof of this fact is due to X.~Wang \cite{WangMass} and Chru\'sciel and Herzlich \cite{ChruscielHerzlich}. Whenever it is non-negative, $m\definedas\sqrt{-\vert\mathbf{P}\vert_{\eta}^{2}}$ is called the \emph{mass} of $(M,g)$\label{def:AHmass}.

In analogy to the asymptotically Euclidean case, the mass vector defined here plays a central role in the \emph{Positive Mass Conjecture}. This conjecture suggests that the mass vector of a geodesically complete asymptotically hyperbolic manifold $(M,g)$ satisfying the (hyperbolic) \emph{dominant energy condition} \label{PMT}
\begin{align*}\labeleq{eq:DEC}
 \scal _g + n(n-1) \geq 0 
\end{align*} 
is timelike future directed, meaning that it satisfies the inequalities
\begin{align*}\labeleq{eq:hypPMT}
 \left(\int_{\bS^{n-1}} \tr _{\sigma} \mathbf{m}\, d\mu_\sigma \right)^2 &>
\sum_{i=1}^n \left(\int_{\bS^{n-1}} \mathring{x}^i \tr _{\sigma} \mathbf{m}\, d\mu_\sigma \right)^2,\\
\int_{\bS^{n-1}} \tr _{\sigma} \mathbf{m}\ d\mu_\sigma &> 0,
\end{align*}
unless $(M,g)$ is isometric to the hyperbolic space $(\bH^{n},b)$. The reader is referred to \cite{ChruscielHerzlich, WangMass, AnderssonCaiGalloway} for proofs of this conjecture under additional assumptions.

\medskip

In the asymptotically Euclidean setting described, one definition of the (total) center of mass of an asymptotically Euclidean Riemannian manifold relies on the existence and uniqueness of a constant mean curvature foliation of its asymptotic end, see Section \ref{subsec:AECMC}. For asymptotically hyperbolic manifolds as discussed here, such CMC-foliations have been studied by Rigger \cite{Rigger}, and Neves and Tian \cite{NevesTian1,NevesTian2}. In particular, the following result was proven in \cite[Theorem 2.2]{NevesTian2}:

\begin{theorem}[Neves-Tian]\label{thFoliation}
Let $(M,g)$ be a 3-dimensional $l$-asymptotically hyperbolic manifold with $l\geq3$. Assume that there is a diffeomorphism $\varphi$ as in Definition \ref{defWang} such that the mass-aspect function $\tr_{\sigma} \mathbf{m}$ is strictly positive and that
\begin{align*}\labeleq{eqBalanced}
 \int_{\bS^{2}} \mathring{x}^i \tr_{\sigma} \mathbf{m}\, d\mu_\sigma = 0 ,\ i = 1,2,3 .
\end{align*}
Then, outside a compact set, $M$ admits a foliation by stable CMC-spheres. 

The foliation is unique among all foliations for which there exists a constant $C$ such that
\begin{align*}\labeleq{InnerOuterRadii}
 \overline{r}_\Sigma - \underline{r}_\Sigma \leq C
\end{align*}
holds for all leaves $\Sigma$ of the foliation, where $\overline{r}_\Sigma \definedas \inf \{r>0\ \vert\ \varphi(\Sigma) \subset B_r (0) \}$ and $\underline{r}_\Sigma \definedas \sup \{r>0 \ \vert\ \varphi(\Sigma)\subset \bH^{3}\setminus B_r (0)\}$ are the outer and inner radius of the $\varphi$-image of a leaf $\Sigma$, respectively, and $B_r(0)$ denotes the geodesic ball of radius $r$ around the origin of $\bH^3$.
\end{theorem}

\begin{definition}
An asymptotically hyperbolic coordinate chart $\varphi$ satisfying \eqref{eqBalanced} is called 
\emph{balanced}. 
\end{definition}

In fact, similar again to the asymptotically Euclidean setting described in Section~\ref{sect:AE}, Neves and Tian \cite[Theorem 2.2]{NevesTian2} provide further estimates for the leaves $\Sigma$ of the constructed foliation. More specifically, the $\varphi$-image of each leaf $\Sigma$ with sufficiently large inner radius and with surface area
$\vert\Sigma\vert_{g} = :4 \pi \sh^2 \hat r$ can be written as a graph over the coordinate sphere centered at the origin with radius $\hat r$:
\begin{align*}\labeleq{eqGraph0}
 \varphi(\Sigma) = \{(\hat r + f(\mathring{x}) , \mathring{x})\ \vert\ \mathring{x} \in \bS^2\},
\end{align*}
where $f: \bS^2 \to \bR$ satisfies
\begin{align*}\labeleq{eqGraph}
 \sup_{\bS^2} \vert f\vert \leq C e^{- \underline{r}} ,\quad \Vert f\Vert_{C^2 (\bS ^2)} \leq C,
\end{align*}
for some constant $C>0$ independent of $\Sigma$. Furthermore, the leaves of the foliation become more and more 'round' as we approach infinity in the sense that there exists a constant $C>0$ independent of $\Sigma$ such that 
\begin{align*}\labeleq{eq:AHroundness}
 \int_{\Sigma} \left\vert\left(\varphi^{*}\partial _r \right)^T\right\vert^2 \,d\mu_g \leq C e^{-2 \underline{r}}, \quad \int_{\Sigma} \vert\mathring A\vert^2 \,d\mu_g \leq C e^{-4 \underline{r}},
\end{align*}
where $\left(\varphi^{*}\partial _r\right) ^T$ is the tangential part of $\varphi^{*}\partial_r$ with respect to $g$, $\mathring{A}$ is the trace free part of the second fundamental form of $\Sigma$ in $(M,g)$, and $d\mu_g$ is the measure induced on $\Sigma$ by $g$. Note that in view of \eqref{InnerOuterRadii}, the above estimates continue to hold when $\underline{r}$ is replaced by $\hat r$, with possibly an adapted constant.

As described in Section \ref{sect:AE}, the asymptotic CMC-foliation of an asymptotically Euclidean Riemannian manifold can be interpreted as an abstract center of mass. Under suitable decay conditions, this abstract center can be coordinatized by taking the limit to infinity of the Euclidean coordinate centers of the leaves of the CMC-foliation, see again Section \ref{sect:AE}, and one might wonder if this is also true for the 'centers' of the foliation constructed by Neves and Tian. As we will see below, this is indeed the case for a suitable definition of the 'coordinate centers' of the leaves.

Of course, the Euclidean coordinate center $\vec{z}_{S}$ of a surface $S\subset\bR^{3}$ is given by 
\begin{align*}\labeleq{eq:barycenter}
\vec{z}_{S}\definedas\frac{1}{\vert S\vert_{\delta}}\int_{S}\vec{x}\,d\mu_{\delta},
\end{align*}
the affine barycenter of the points on $S$. However, unlike Euclidean space, hyperbolic space is not affine, so first of all we need to clarify how the hyperbolic coordinate center of a surface $S\subset\bH^{3}$ is defined. For this purpose, it is natural to exploit the isometric embedding $I:(\bH^3,b) \hookrightarrow (\bR^{3,1}, \eta)$ of the hyperbolic space into the Minkowski spacetime $(\bR^{3,1}, \eta)$ given in canonical coordinates by
\begin{equation*}
I:(r,\mathring{x}^1, \mathring{x}^2, \mathring{x}^3) \mapsto (\ch r, \mathring{x}^1 \sh r, \mathring{x}^2 \sh r, \mathring{x}^3\sh r),
\end{equation*}
where $\mathring{x}^i$ denotes again the restriction to $\bS^2$ of the Cartesian coordinate $x^i$ in $\bR^3$. We can therefore introduce 'Minkowskian' coordinates $(X^{\alpha})=(X^0,X^1,X^2,X^3)$ on $\bR^{3,1}$ so that, on $I(\bH^3)$, we have $X^{\alpha} = I^{\alpha}(r,\mathring{x}^1, \mathring{x}^2, \mathring{x}^3)$ for $\alpha = 0,1,2,3$ and attempt to view the hyperbolic center of a surface $S\subset\bH^{3}$ as the vector $\mathbf{C} = (C^0, C^1, C^2, C^3)\in\bR^{3,1}$ with components
\begin{align*}\labeleq{eqCenterMink}
 C^\alpha \definedas \frac{1}{\vert S\vert_{b}}\int_{S} I^\alpha\, d\mu_{b} = \frac{1}{\vert I(S)\vert_{I_* b}} \int_{I(S)} X^{\alpha} d\mu_{\eta}, \quad \alpha = 0,1,2,3.
\end{align*}
This Minkowskian 'definition', however, is too naive: for the centered CMC-spheres of hyperbolic space, $S_{R}=\{r=R\}\subset\bH^3$, formula \eqref{eqCenterMink} gives $\mathbf{C}_{R} = (\ch R, 0, 0, 0)$. This neither converges to a 'center of mass' of the manifold as $R\to\infty$, nor does it correspond to a 'point' in the $I$-image of hyperbolic space. It thus seems more natural to define the \emph{hyperbolic center $\mathbf{z}$ of a surface $S \subset \bH^{3}$} as 
\begin{align*}\labeleq{eqCenterHyperboloid}
\mathbf{z}\definedas I^{-1}(\mathbf{Z}) \definedas I^{-1}\left(\frac{\mathbf{C}}{\sqrt{-\vert \mathbf{C}\vert_\eta}}\right),
\end{align*}
where $\mathbf{C}$ is as above. Defined in this way, the hyperbolic center of $S$ can be thought of as a point $\mathbf{z}\in\bH^3$, or, equivalently, 
as a future timelike unit vector $\mathbf{Z}\in\bR^{3,1}$. 

\begin{remark}\label{remGeneralcoordcenter}
In fact, our definition of the hyperbolic center of a surface $S\subset\bH^{3}$ is a special case of a more general idea which is closely related to affine geometry and special relativity and also ties in with the definition of quasi-local center of mass by Chen, Wang, and Yau, see Section \ref{sect:AE}. To see this, assume we are given a 'model space' Riemannian $n$-manifold $(\mathbb{M}^{n},h)$ that can be isometrically embedded into the Minkowski spacetime $(\bR^{n,1},\eta)$ via an embedding $I:(\mathbb{M}^{n},h)\hookrightarrow(\bR^{n,1},\eta)$ that is \emph{affinely bijective}, meaning that for each $\mathbf{X}\neq\mathbf{0}$ in the convex cone $\cone(I(\mathbb{M}^{n}))$ spanned by $I(\mathbb{M}^{n})$ with apex at the origin, one finds a unique $\mathbf{z}\in\mathbb{M}^{n}$ such that $I(\mathbf{z})\,\Vert\,\mathbf{X}$. Then, the corresponding center of any closed hypersurface $S^{n-1}\subset\mathbb{M}^{n}$ can be defined as the unique $\mathbf{z}\in\mathbb{M}^{n}$ such that $I(\mathbf{z})\,\Vert\,\mathbf{C}$, where the components of $\mathbf{C}$ are given by
\begin{align*}\labeleq{eq:generalCoM}
C^{\alpha}\definedas \int_{S} I^\alpha\, d\mu_{h}= \int_{I(S)} X^\alpha\, d\mu_{\eta}, \quad \alpha = 0,1,2,\dots,n.
\end{align*}
Indeed, $\mathbf{C}\in\cone(I(\mathbb{M}^{n}))$ as $I(S)\subset I(\mathbb{M}^{n})\subset\cone(I(\mathbb{M}^{n}))$ and because integration is an affine operation. 

In particular, if the model space is $(\bH^{3},b)$ and the embedding $I$ is as above, then $\cone(I(\mathbb{M}^{n}))$ is the (open) light cone and $I$ is thus clearly affinely bijective. For any closed surface $S\subset\bH^{3}$, the unique $\mathbf{z}\in\bH^{3}$ with $I(\mathbf{z})\,\Vert\,\mathbf{C}$ is precisely found by \eqref{eqCenterHyperboloid}, as the embedding $I$ maps $(\bH^{3},b)$ onto the unit hyperboloid in Minkowski space. Similarly, if the model space is $(\bR^{3},\delta)$ and $I:{x}\mapsto(\tau,{x})$ for any fixed $\tau\neq0$, $\cone(I(\mathbb{M}^{n}))$ is one of the half-spaces $\{t\lessgtr0\}$, depending on the sign of $\tau$, and $I$ again is of course affinely bijective. We find
\begin{align*}\labeleq{eq:C-components}
\mathbf{C}=\left(\tau\,{\vert S\vert_{\delta}},\int_{S}{x}\,\,d\mu_{\delta}\right)
\end{align*}
so that $\mathbf{C}\,\Vert\,I(\vec{z})$ precisely when $\vec{z}$ is the Euclidean center of $S$.
\end{remark}

\begin{remark}\label{rem:transfo}
This definition of the center of a surface transforms appropriately under \emph{time-preserving} Lorentz transformations of the ambient Minkowski space because they are affine. Intuitively, one can understand this definition as saying that the center of $S$ is the $I$-preimage of the Minkowskian center $\mathbf{C}$ of $I(S)$ if $\mathbf{C}\in I(\mathbb{M}^{n})$, and the unique $I$-preimage corresponding to the line through the origin in direction $\mathbf{C}$ in general. This generalizes the affine concept of a Euclidean center of a surface to a not necessarily affine ambient 'model space' which can be isometrically embedded into Minkowski space\footnote{In fact, any flat $\bR^{l,n+1-l}$ could also be chosen as a reference space.} by exploiting the extrinsic affine structure induced by the isometric embedding.
\end{remark}

Equipped with this definition of a hyperbolic center of a surface $S\subset\bH^{3}$, we can now go back to defining the hyperbolic coordinate center corresponding to the asymptotic CMC-foliation of an asymptotically hyperbolic Riemannian manifold. As in the asymptotically Euclidean setting, we say that the \emph{hyperbolic coordinate center} of an $l$-asymptotically hyperbolic Riemannian manifold $(M^{3},g)$ as in Theorem \ref{thFoliation} and with respect to the chart $\varphi$ is the limit of the hyperbolic centers $\mathbf{z}_{\varphi(\Sigma_{\hat{r}})}= I^{-1}(\mathbf{Z}_{\varphi(\Sigma_{\hat{r}})})$ of the ($\varphi$-images of the) leaves $\Sigma_{\hat{r}}$ of the CMC-foliation as the area radius $\hat{r}\to\infty$, whenever this limit exists. Consistently, for the centered CMC-spheres $S_{R}=\{r=R\}\subset\bH^3$, we have $\mathbf{Z}_{S_{R}}= \mathbf{N}$, where $\mathbf{N}=(1,0,0,0)$, and thus $\mathbf{z}_{S_{R}}=\mathbf{0}$. Thus the centers of the leaves of this foliation (trivially) converge to $\mathbf{0}$, the coordinate origin of hyperbolic space. Using the graphical representation \eqref{eqGraph0} together with the estimates \eqref{eqGraph}, we will reach the same conclusion for the foliation of Neves and Tian: not only does the hyperbolic coordinate center always converge under the assumed fall-off conditions, it indeed always lies at the chosen origin of hyperbolic space:

\begin{proposition}\label{propConvergence}
Under the conditions of Theorem \ref{thFoliation}, the hyperbolic centers of the leaves of the constant mean curvature foliation constructed in Theorem \ref{thFoliation} converge to $\mathbf{0}$. In other words, for $l \geq 3$, the hyperbolic coordinate center of an $l$-asymptotically hyperbolic Riemannian manifold with strictly positive mass aspect function and balanced coordinates lies at the center of the hyperbolic coordinates.
\end{proposition}

\proof
Let $\Sigma$ be a leaf of the CMC-foliation in Theorem \ref{thFoliation} with surface area $\vert\Sigma\vert_{g} \asdefined 4 \pi \sh^2 \hat r$ and $\hat r$ sufficiently large. Then $\varphi(\Sigma)$ can be written as a graph over the coordinate sphere of radius $\hat r$ as in \eqref{eqGraph0}. As a consequence, we may write $X^0\vert_{\varphi(\Sigma)} = \ch(\hat r + f)$, and $X^i \vert_{\varphi(\Sigma)} = \sh (\hat r + f)\, \mathring{x}^i$, $i=1,2,3$, where $f : \bS^2 \to \bR$ satisfies \eqref{eqGraph}. Using the estimates \eqref{eqGraph} -- which are uniform in $\hat{r}$ -- and a Taylor expansion near $f=0$ we see that on $\varphi(\Sigma)$
\begin{align*}\labeleq{eq:Taylor}
 g^\phi\vert_{\varphi(\Sigma)} &= \sh^2 \hat r \left[\sigma + \mathcal{O}(e^{- \hat{r}})\right],\\
 d\mu_{g^{\phi}} &= \sh^2 \hat r \left(1 + \mathcal{O}(e^{- \hat{r}})\right) d\mu_\sigma.
\end{align*}
Note also that $\sh(\hat{r} + f) = \sh \hat r \left(1 + \mathcal{O}(e^{-\hat r})\right)$. Using \eqref{eqCenterMink}, we compute
\begin{align*}\labeleq{eq:C-components-AH}
C^0 & = \frac{1}{4\pi } \int_{\bS^2} \ch \hat r \left(1 + \mathcal{O}(e^{- \hat{r}})\right)d\mu_\sigma \phantom{\mathring{x}^{i}}= \ch \hat r + \mathcal{O}(1),\\
 C^i & = \frac{1}{4\pi } \int_{\bS^2} \mathring{x}^i \sh \hat r \left(1 + \mathcal{O}(e^{- \hat{r}})\right)d\mu_\sigma = \mathcal{O}(1),
 \end{align*}
as $\hat r\to\infty$, for $i=1,2,3$. Then $\sqrt{-\vert\mathbf{C}\vert_\eta} = \ch \hat r + \mathcal{O}(1)$, and thus the components of $\mathbf{Z}$ satisfy $Z^0= 1+\mathcal{O}(e^{-\hat r})$, $Z^i = \mathcal{O}(e^{-\hat r})$, $i=1,2,3$ so that the hyperbolic centers $\mathbf{z}=\mathcal{O}(e^{-\hat r})$ of the leaves tend to the hyperbolic coordinate center $\mathbf{0}$.
\qed

\begin{remark}\label{remarkBoost} 
Proposition \ref{propConvergence} shows that the condition that the asymptotic coordinate chart $\varphi$ be balanced forces the hyperbolic coordinate center of mass to lie at the center of the coordinate system, which is possibly not too surprising. However, both the existence and uniqueness of the asymptotic CMC-foliation and the definition of its hyperbolic coordinate center naturally extend to asymptotic charts $\widetilde{\varphi}$ that satisfy all assumptions of Theorem \ref{thFoliation} but that are not necessarily balanced, see Theorem \ref{thCoMUnbalanced}.
\end{remark}

\begin{remark}
At this point it is also natural to ask whether the unique constant mean curvature foliation still exists if we merely assume that the mass of $(M,g)$ is 'positive', meaning that the mass vector $\mathbf{P}\in \bR^{3,1}$ is timelike future directed. This seems to be a non-trivial question. Indeed, the assumption that the mass aspect function is strictly positive plays a crucial role in \cite{NevesTian2}. In particular, it is used by Neves and Tian \cite[Section 8]{NevesTian2} to ensure that the \emph{normalized Jacobi operator} is positive definite.
\end{remark}
 \label{Riemannian}

\subsection{General situation}\label{sect:CoM}\label{sect:geometric}
% !TEX root = CenterofMassAH.tex
As discussed in Section \ref{sect:balanced}, Theorem \ref{thFoliation} of Neves and Tian guarantees that the \emph{abstract center of mass} is well-defined for manifolds which are asymptotically hyperbolic in the sense of Definition \ref{defWang} with respect to a \emph{balanced} chart at infinity and which have strictly positive mass aspect function. We showed in Proposition \ref{propConvergence} that in this case the hyperbolic centers of the leaves of the constant mean curvature foliation do converge when the area of the leaves tends to infinity. In fact, both results hold true without assuming that the chart at infinity is balanced. More precisely, we have: 

\begin{theorem}\label{thCoMUnbalanced}
Let $(M,g)$ be a 3-dimensional $l$-asymptotically hyperbolic manifold in the sense of Definition \ref{defWang} for $l\geq 3$, with a chart at infinity denoted by $\phi$. Assume that the mass aspect function $\tr_{\sigma} \mathbf{m}$ is strictly positive. Then outside of a compact set, $M$ admits a foliation by stable CMC-spheres. This foliation is unique among all foliations which have Property \eqref{InnerOuterRadii}.

Furthermore, let $\mathbf{P}$ be the mass vector of $(M,g)$ and let $\mathbf{z}_{\hat r} \in \bH^3$ be the hyperbolic center of the leaf $\Sigma_{\hat r}$ with surface area $4\pi \sh^2 \hat r$ as in \eqref{eqCenterMink}. Then 
 \begin{align*}\labeleq{eq:unbalanced-center}
  \mathbf{z}_{\hat r} \longrightarrow \mathbf{p}\definedas I^{-1}\left(\frac{\mathbf{P}}{\sqrt{-\vert \mathbf{P}\vert _{\eta}^2}}\right) \quad \text{as} \quad \hat r \rightarrow \infty.
 \end{align*}
\end{theorem}

\proof
Let $(M,g)$ and $\phi$ be as in the statement of the theorem. Since the mass aspect function is positive, it follows from 
the Cauchy-Schwarz inequality -- with respect to the measure $\tr_{\sigma} \mathbf{m}\, d\mu_{\sigma}$ -- that the mass vector $\mathbf{P}$ is timelike future directed in $\bR ^{3,1}$, see also Cortier \cite[Remark 2.10]{Cortier}.

Since the group of isometries $SO_{0}(3,1)$ acts transitively on $(\bH^3,b)$, there exists $A_{\mathbf{P}}\in SO_{0}(3,1)$ such that
\begin{align*}\labeleq{eq:A_P}
A _{\mathbf{P}}:\mathbf{P} \mapsto \sqrt{-\vert \mathbf{P}\vert _{\eta}^2}\ \mathbf{N}.
\end{align*}
Then $A_{\mathbf{P}}$ is a (linear) isometry of $\bR^{3,1}$ whereas its restriction $a\definedas A_{\mathbf{P}}\circ I$ to the hyperbolic space $I:\bH^3 \hookrightarrow \bR^{3,1}$ is non-linear. The mass vector with respect to the chart $a\circ \phi$ at infinity is $A_{\mathbf{P}}(\mathbf{P})$. In other words, the chart $a\circ\phi$ is {balanced}. Let us also note that the condition \eqref{InnerOuterRadii} is preserved under the hyperbolic isometry $a$. Furthermore, the mass-aspect function of $g^{a\circ\phi} := (a\circ \phi)_* g$ is positive. This follows from an observation by Wang \cite[p.~291]{WangMass}: the round metric transforms conformally under the action of $a\vert_{\bS^{2}}$ seen as a conformal diffeomorphism of $\bS^2$, and so does the mass-aspect tensor with a different power of the (positive) conformal factor $u:\bS^{2}\to\bR$:
\begin{align*}\labeleq{eq:a}
 a \cdot \sigma = a_*\sigma = u^{-2}\, \sigma\ ,\ a \cdot \mathbf m = u\; a_*\mathbf m\ ,\ a \cdot \tr_{\sigma} \mathbf m = u^3\, \tr_{\sigma} \mathbf m \circ a^{-1},
\end{align*}
where the ``$\cdot$'' denote the respective actions of conformal diffeomorphisms, see Cortier, Dahl, and Gicquaud \cite[Section 3]{CDG} for details. Consequently, Theorem \ref{thFoliation} applies with respect to the chart $a\circ \phi$ to produce a foliation by constant mean curvature spheres. By Proposition \ref{propConvergence}, these spheres are asymptotically hyperbolically centered at $\mathbf{0} \in \bH^3$. 

However, the CMC-foliation is coordinate independent and thus exists without reference to the balanced coordinates $a\circ\varphi$. From the discussions in Remarks \ref{remGeneralcoordcenter} and \ref{rem:transfo}, the hyperbolic coordinate center with respect to the chart $\phi$ of the CMC-foliation is at $\mathbf{z}=I^{-1}(A_{\mathbf{P}}^{-1}(\mathbf{N}))$ with
\begin{align*}\labeleq{eq:CMC-hyp-center}
 A_{\mathbf{P}}^{-1}\mathbf{N} &= \tfrac{1}{\sqrt{-\vert \mathbf{P}\vert _{\eta}^2}}A_{\mathbf{P}}^{-1}\left(\sqrt{-\vert \mathbf{P}\vert _{\eta}^2}\ \mathbf{N}\right) = \tfrac{\mathbf{P}}{\sqrt{-\vert \mathbf{P}\vert _{\eta}^2}}=I(\mathbf{p})
\end{align*}
so that $\mathbf{z}=\mathbf{p}$, regardless of the choice of the element $A_{\mathbf{P}}$ made earlier.
\qed

This motivates us to adopt the following general definition of center of mass.

\begin{definition}\label{defCoMRiem}
 Let $(M,g)$ be an $n$-dimensional asymptotically hyperbolic manifold in the sense of Definition \ref{defWang} with respect to a chart $\varphi$. Assume that it has \emph{positive mass}, in the sense that the mass vector $\mathbf{P}$ of $(M,g)$ is timelike future directed. Then the \emph{hyperbolic coordinate center of mass} $\mathbf{z}$ with respect to the chart $\phi$ is the given by
 \begin{align*}\labeleq{eq:hypCoM}
  \mathbf{z}\definedas I^{-1}\left(\frac{\mathbf{P}}{\sqrt{-\vert \mathbf{P}\vert _{\eta}^2}}\right) \in\bH^{n},
 \end{align*}
 where $I:(\bH^{n},b)\hookrightarrow(\bR^{n,1},\eta)$ is the canonical isometric embedding.
\end{definition}

\begin{remark}
 This definition interprets the center of mass as a point in hyperbolic space. Note also that it transforms appropriately under the action of the ($I$-induced) intrinsic isometry group $SO_0(n,1)$, which is consistent with the discussion in Remarks \ref{remGeneralcoordcenter} and \ref{rem:transfo}. In particular, these properties are reminiscent of the properties of the asymptotically Euclidean center of mass, see Section \ref{subsec:AEprop}.
\end{remark}

\section{Hamiltonian charges of asymptotically hyperbolic\\ initial data sets}
% !TEX root = CenterofMassAH.tex
In Section \ref{sect:AAdS}, we have defined the hyperbolic coordinate center of an asymptotically hyperbolic Riemannian $n$-manifold of positive mass, see Definition \ref{defCoMRiem}. We have argued that, for $n=3$, the center we define coincides with the limit of the hyperbolic centers of the leaves of the CMC-foliation in the asymptotic end, just as in the asymptotically Euclidean setting. Also, for general $n\geq3$, we noted that the hyperbolic coordinate center transforms adequately under the symmetry group of the hyperbolic space, in analogy to the transformation of the Euclidean coordinate center of mass under Euclidean motions. In this section, we will investigate how the hyperbolic coordinate center relates to the Hamiltonian charges derived for asymptotically hyperbolic initial data sets of asymptotically anti-de Sitter spacetimes.

For this, recall that the hyperbolic space $(\bH^n,b)\cong (\bR\times \bS^{n-1}, dr^2 + \sh^2 r\, \sigma)$ naturally arises as the time-symmetric spacelike hypersurface $\{t = 0\}$ in the anti-de Sitter spacetime 
\begin{align*}\labeleq{eqAdS}
\left(\mathrm{AdS} ^{n+1}, \mathfrak{g}_{\mathrm{AdS}}\right) = (\bR \times \bH^{n},  - \ch^2 r\, dt^2 + b).
\end{align*} 
It may therefore be viewed as the initial data set $(\bH^n, b, 0)$ for $\left(\mathrm{AdS} ^{n+1}, \mathfrak{g}_{\mathrm{AdS}}\right)$. In general, we will define an asymptotically hyperbolic (asymptotically anti-de Sitter) initial data set as follows.\\

\begin{definition}\label{defAdSData}
Let $(M,g)$ be a Riemannian $n$-manifold and let $k$ be a symmetric $2$-tensor field on $M$. We say that $(M,g,k)$ is an \emph{asymptotically hyperbolic initial data set of order $\tau>0$} if there exist a compact set $K \subset M$, a radius $R>0$, and a diffeomorphism $\phi: M\setminus K \to \bH^n \setminus \overline{B}_R$ such that $\gamma \definedas \phi_*g - b = \mathcal{O}_2 (e^{-\tau r})$ and $\kappa \definedas \phi_* k =\mathcal{O}_1 (e^{-\tau r})$. 
\end{definition}

\begin{remark}
This generalizes Definition \ref{defWang} in the sense that if $(M,g)$ is an asymptotically hyperbolic Riemannian $n$-manifold according to Definition \ref{defWang} then the initial data set $(M,g,0)$ is an asymptotically hyperbolic initial data set of order $\tau=n$ according to Definition \ref{defAdSData}.
\end{remark}

Given a Killing vector field $X$ of a background spacetime, there is a standard procedure in General Relativity which allows to define a \emph{Hamiltonian charge} associated with the flow along $X$, see also Section \ref{sect:AE} for the case of the Minkowski background spacetime. This approach was used by Chru\'sciel and Nagy \cite{ChruscielNagy} to define what they call the 'global charges' for asymptotically hyperbolic initial data sets as in Definition \ref{defAdSData}. We briefly recall this construction following Michel \cite{MichelMass}. Given a Killing vector field $X$ of the anti-de Sitter spacetime \eqref{eqAdS}, the associated \emph{Killing initial data} (or \emph{KID}) on the slice $\{t=0\}$ is a pair $(V,Y)$, where the 
\emph{lapse} function $V$ and the \emph{shift} vector $Y$ are uniquely determined by the formula
\begin{align*}\labeleq{eq:lapseshift}
 X\vert _{\{t=0\}} &\asdefined V n + Y.
\end{align*}
When the background spacetime is $\mathrm{AdS}^{n+1}$, $n$ is the future directed timelike unit normal to the slice $\{t=0\}\cong (\bH^n,b,0)$ in $\mathrm{AdS}^{n+1}$, $V: \bH^n \rightarrow \bR$, and $Y$ is a tangent vector field along $\bH^{n}$. Now let $(M,g,k)$, $\phi$, $\gamma$, and $\kappa$ be as in Definition \ref{defAdSData}. Then to every KID $(V,Y)$, one may formally assign a \emph{global charge} 
\begin{align*}\labeleq{eq:Q}
 \mathbb{Q}_{(V,Y)} (\gamma,\kappa) &\definedas \lim_{R \to \infty} \int_{\{r=R\}} \bU_{(V,Y)}(\gamma,\kappa)(\nu_{b})\,d\mu_b,
\end{align*}
where the 1-form $\bU_{(V,Y)}(\gamma,\kappa)$ is given by the expression
\begin{align*}\labeleq{eqChargeDensity}
 \bU_{(V,Y)}(\gamma,\kappa) &\definedas V \left(\divg_b \gamma - d (\tr_b \gamma)\right) - \iota _{\nabla_b V} \gamma + \tr_b \gamma\ dV
 + 2\left(\iota _{Y}\kappa - \left(\tr_b \kappa\right) Y^\flat\right)
\end{align*}
and $\nu_{b}$ being the outward unit normal of $\{r=R\}$ with respect to $b$. Under suitable decay conditions on the lapse $V$ and the shift $Y$, the formally defined global charges converge; these conditions will be recalled in Proposition \ref{propWellDefined}.

To make this construction more explicit, we need to describe KIDs of the anti-de Sitter spacetime. For convenience, we will switch to the Poincar\'e ball model of hyperbolic space for a while, where $(\bH^n,b)$ is seen as the unit ball 
\begin{align*}
\{y\in (\bR^n, \delta) : \vert y\vert  < 1\}
\end{align*}
and the metric is given by 
\begin{align*}\labeleq{eq:ballmetric}
 b &= \left(\frac{2}{1-\vert y\vert ^2}\right)^2 \delta .
\end{align*}
Note that the radial coordinate $r$ corresponding to the polar coordinates used to describe the hyperbolic space before is related to the ball coordinates $(y^{i})$ via the formula 
\begin{align*}
r &= \mathrm{arcch}\, \frac{1+\vert y\vert ^2}{1-\vert y\vert ^2}.
\end{align*}
Hence, the hyperbolic metric in these coordinates reads 
\begin{align*}\labeleq{eqAdSBall}
\mathfrak{g}_{\mathrm{AdS}} = -\left(\frac{1+\vert y\vert ^2}{1-\vert y\vert ^2}\right)^2 dt^2 + b.
\end{align*}
The algebra of Killing vector fields of $(\mathrm{AdS}^{n+1},\mathfrak{g}_{\mathrm{AdS}})$ is given by 
\begin{align*}\labeleq{KVFs}
\begin{split}
X_{(0)} & \definedas\frac{\partial}{\partial t}\ ,\\
X^+_{(i)} & \definedas\frac{2y^i \cos t}{1+\vert y\vert ^2}\frac{\partial}{\partial t} + \left( \frac{1+\vert y\vert ^2}{2} \delta^j_i -y^i y^j\right)\sin t\,\frac{\partial}{\partial y^j}\ , \hspace{0.5cm} i=1,\ldots, n, \\
X^-_{(i)} & \definedas
 \frac{2y^i \sin t}{1+\vert y\vert ^2}\frac{\partial}{\partial t} - \left( \frac{1+\vert y\vert ^2}{2} \delta^j_i -y^i y^j\right)\cos t\,\frac{\partial}{\partial y^j} \ , \hspace{0.5cm}  i=1,\ldots, n,  \\
 X_{(i)(j)} & \definedas y^i \frac{\partial}{\partial y^j} - y^j \frac{\partial}{\partial y^i}, \hspace{0.5cm}  1\leq i < j \leq n.
\end{split}
\end{align*}
The future timelike unit normal vector field along the slice $\{t=0\}$ is $\nu = \frac{1-\vert y\vert ^2}{1+\vert y\vert ^2} \frac{\partial}{\partial t}$. From this, one obtains the following lapse functions and shift vectors:
\begin{figure}[h]
\begin{center}
\begin{tabular}{|c|| c| c| }
Killing V.F. & Lapse & Shift \\\hline\hline
&&\\[-1ex]
$X_{(0)}$ & $ V_{(0)}=\frac{1+\vert y\vert ^2}{1-\vert y\vert ^2} $ & $0$ \\[1ex]\hline
&&\\[-1ex]
$ X^+_{(i)}   $ & $V_{(i)} =  \frac{2y^i}{1-\vert y\vert ^2} $ & $0$ \\[1ex]\hline
&&\\[-1ex]
$ X^-_{(i)} $ & $ 0 $ & $ C_{(i)} \definedas -\frac{1+\vert y\vert ^2}{2} \frac{\partial}{\partial y^i} + y^i y^j \frac{\partial}{\partial y^j}$ \\[1ex]\hline
&&\\[-1ex]
$ X_{(i)(j)}$ & $0$ & $\Omega_{(i)(j)} \definedas y^i \frac{\partial}{\partial y^j} - y^j \frac{\partial}{\partial y^i}$\\[1ex]\hline
\end{tabular}
\end{center}
\caption{Killing initial data (KIDs) for the $\{t=0\}$-slice of the AdS spacetime corresponding to the Killing vector fields \eqref{KVFs} in the ball model of hyperbolic space.}\end{figure}

In \cite[Section 3]{ChruscielMaertenTod}, Chru\'sciel, Maerten and Tod gave a physical interpretation of the respective global (Hamiltonian) charges, namely
\begin{itemize}
\item The \emph{energy-momentum vector} was defined as
\begin{align*}\labeleq{eqMassVector}
\left(\mathbb{Q}_{(V_{(0)},0)}, \mathbb{Q}_{(V_{(1)},0)}, \ldots, \mathbb{Q}_{(V_{(n)},0)} \right).
\end{align*}
\item The vector 
\begin{align*}\labeleq{eqCoM}
\left(\mathbb{Q}_{(0, C_{(1)})}, \ldots, \mathbb{Q}_{(0,C_{(n)})}\right)
\end{align*}
 was associated with the \emph{center of mass}.
\item The global charges 
\begin{align*}\labeleq{eqAngularMomentum}
\mathbb{Q}_{(0,\Omega_{(i)(j)})},\;  1\leq i < j \leq n,
\end{align*}
were identified with components of the \emph{angular momentum}.
\end{itemize}

Changing back to the coordinates $(r,\mathring{x}^1, \ldots, \mathring{x}^n)$, we obtain
\begin{align*}\labeleq{eq:vballKIDs}
\begin{split}
V_{(0)} &= \ch r\\
V_{(i)} &= \mathring{x}^i \sh r
\end{split}
\end{align*}
for $i=1, \ldots , n$. If $(M,g)$ is an asymptotically hyperbolic manifold in the sense of Definition \ref{defWang}, it is straightforward to check that the energy-momentum vector \eqref{eqMassVector} of the initial data set $(M,g,0)$ coincides with what we called the mass vector \eqref{eqMassVectorRiem} of $(M,g)$, up to a dimensional multiple. By analogy, the vector \eqref{eqMassVector} will be referred to as the \emph{mass vector} of the initial data set $(M,g,k)$ in what follows. Note that this vector transforms equivariantly under the action of $SO_0(n,1)$, see also Section~\ref{sect:balanced}.

At the same time, it is clear that the vector \eqref{eqCoM} is trivial in the case of the initial data set $(M,g,0)$, whereas the coordinate center of the CMC-foliation\footnote{provided that it exists} does not need to be located at the coordinate origin $\mathbf{0}$. In other words, the hyperbolic coordinate center and the Hamiltonian notion of the center of mass do not coincide.

Generalizing Definition \ref{defCoMRiem}, we thus suggest the following definition of the center of mass for asymptotically hyperbolic initial data sets of positive mass. Further motivation will be given by Theorem \ref{thEvolution} which asserts that the hyperbolic center of mass evolves appropriately under the Einstein evolution equations. In particular, we will see that the so-defined center of mass evolves in the direction of the 'linear momentum' which we define via \eqref{eqCoM}, see Section \ref{Evolution}.

\begin{definition}\label{defCoMcharges}
 Let $(M,g,k)$ be an asymptotically hyperbolic initial data of order $\tau$ as in Definition \ref{defAdSData} with respect to a chart $\phi$ near infinity. Assume that the mass vector
\begin{align*}\labeleq{eq:massvector}
\mathbf{P} &\definedas \left(\mathbb{Q}_{(V_{(0)},0)}, \mathbb{Q}_{(V_{(1)},0)}, \ldots, \mathbb{Q}_{(V_{(n)},0)} \right)\in\bR^{n,1}
\end{align*}
is well-defined and timelike future directed in the Minkowski space $(\bR^{n,1},\eta)$. Then the \emph{center of mass} $\mathbf{z}$ with respect to the chart $\phi$ is given by
 \begin{align*}\labeleq{eq:chargesCoM}
  \mathbf{z} \definedas I^{-1}\left(\frac{\mathbf{P}}{\sqrt{-\vert \mathbf{P}\vert _\eta}}\right) \in \bH^n,
 \end{align*}
 where $I:(\bH^{n},b)\hookrightarrow(\bR^{n,1},\eta)$ is the canonical isometric embedding.
\end{definition}

Again, the center of mass is a point in the hyperbolic space and it transforms adequately under the action of the group $SO_0(n,1)$. Note that the positive mass theorem states that the mass vector of a geodesically complete asymptotically hyperbolic initial data set in an asymptotically AdS spacetime satisfying the dominant energy condition $\mu\geq\vert J\vert$ must be future timelike, unless the spacetime is precisely the AdS spacetime, compare p.~\pageref{PMT} and \cite{Maerten}. Here, $\mu$ and $J$ are defined as below.

We conclude the current section by recalling assumptions which guarantee that the quantities \eqref{eqMassVector}-\eqref{eqAngularMomentum} are well-defined. Let the \emph{energy density} $\mu$ and the \emph{momentum density} $J$ be defined through the  \emph{constraint equations} with cosmological constant $\Lambda=-\frac{n(n-1)}{2}$ as 
\begin{align*}
\mu &\definedas \scal_g + n(n-1) -(\tr_g k)^2 + \vert k\vert _g^2 \,,\labeleq{eqConstraints1}\\ 
J &\definedas\divg_g k -d (\tr_g k). \labeleq{eqConstraints2}
\end{align*}
The reader is referred to Michel \cite{MichelMass} for the proof of the following fact.
\begin{prop}\label{propWellDefined}
 Let $(M,g,k)$ be an asymptotically hyperbolic initial data set of order $\tau > \tfrac{n}{2}$ as in Definition \ref{defAdSData} with respect to a chart $\varphi$ near infinity. If
 \begin{align*}\labeleq{eq:chargewelldef}
\langle (\phi_*\mu, \phi_*J), (V,Y) \rangle \in L^1 (d\mu_b),
 \end{align*}
 then the charge $\mathbb{Q}_{(V,Y)}(\gamma,\kappa)$ is well-defined. 
\end{prop}
 
Note that for any KID $(V,Y)$ for anti-de Sitter spacetime we have $V = \mathcal{O}(e^{r})$, $ Y  = \mathcal{O}(e^r)$ as $r\to\infty$. In particular, the charges \eqref{eqMassVector}-\eqref{eqAngularMomentum} will be well-defined provided that 
$\phi_* \mu = \mathcal{O}(e^{-(n+\epsilon)r})$ and $ \phi_* J  = \mathcal{O}(e^{-(n+\epsilon)r})$ for some $\epsilon > 0$.\label{Charges}

\section{Evolution of the hyperbolic center of mass}\label{Evolution}
% !TEX root = CenterofMassAH.tex
\noindent To further justify our definition of center of mass for asymptotically hyperbolic initial data sets, we will now show that it evolves appropriately under the Einstein evolution equations, see also Section \ref{sect:AE}. For simplicity, we will from now on suppress the chart $\varphi$ near infinity, and identify a tensor field $T$ on $M \setminus K$ with its pushforward $\varphi_* T$ on $\bH^n \setminus \overline{B}_R$. In particular, we will write $g, k, \mu, J$ instead of $\phi_*g, \phi_* k$, $\phi_* \mu, \phi_* J$.

\begin{theorem}\label{thEvolution}
Let $(M^{n},g,k)$ be an asymptotically hyperbolic initial data set of order $\tau>\frac{n}{2}$ satisfying the constraint equations \eqref{eqConstraints1}-\eqref{eqConstraints2} with negative cosmological constant $\Lambda\definedas-\frac{n(n-1)}{2}$, $\mu=\mathcal{O}\left(e^{-(n+\varepsilon)r}\right)$, and $J = \mathcal{O}\left(e^{-(n+1+\varepsilon)r}\right)$ for some $\varepsilon > 0$. Assume additionally that $k  = \mathcal{O}\left(e^{-(\tau+1) r }\right)$. Consider a family $(g(t),k(t))$ that evolves starting from $(g(0),k(0))=(g,k)$ according to the Einstein evolution equations with cosmological constant $\Lambda$, lapse $N=V_{(0)}$ and shift $X=0$. Then
\begin{align*}
\left.\frac{d}{dt}\right\vert _{t=0}\mathbb{Q}_{(V_{(0)},0)} &= 0,\\
\left.\frac{d}{dt}\right\vert _{t=0}\mathbb{Q}_{(V_{(i)},0)} &= \mathbb{Q}_{(0,C_{(i)})}, \quad i=1,\ldots, n.
\end{align*}
\end{theorem}

\proof
Note that by Proposition \ref{propWellDefined}, all global charges involved in the formulation of this theorem are well-defined. 
By the Einstein evolution equations\footnote{see for example Bartnik and Isenberg \cite{BartnikIsenberg}.} with cosmological constant $\Lambda=-\frac{n(n-1)}{2}$, we have
\begin{align*}
\dot{g} \definedas \left.\frac{d}{dt}\right\vert _{t=0}g = 2 V_{(0)} k = \mathcal{O}(e^{-\tau r}).
\end{align*}
Using  the momentum constraint \eqref{eqConstraints2}, we compute for $\beta = 0, \ldots, n$ 
\begin{align*}
\begin{split}
& \left.\frac{d}{dt}\right\vert _{t=0}\, \bU_{(V_{(\beta)},0)} (g,k) \\
& \quad = V_{(\beta)} (\divg \dot{g} - d \tr \dot{g}) - \iota_{\nabla V_{(\beta)}} \dot{g} + \tr \dot{g}\, dV_{(\beta)} \\ 
& \quad = 2V_{(\beta)} \left( \divg (V_{(0)}k) - d (V_{(0)}\tr k)\right) -2 V_{(0)} \left(\iota_{\nabla V_{(\beta)}} k - \tr k\, dV_{(\beta)}\right) \\
& \quad = 2V_{(\beta)} \left(\iota_{\nabla V_{(0)}} k + V_{(0)} \divg k - \tr k \, d V_{(0)} - V_{(0)} d (\tr k)\right) \\
& \quad\quad -2 V_{(0)}\left( \iota_{\nabla V_{(\beta)}}k - \tr k\, dV_{(\beta)} \right)\\
& \quad = 2 \iota_{V_{(\beta)} \nabla V_{(0)} - V_{(0)} \nabla V_{(\beta)}} \left(k - (\tr k) b\right) + 2 V_{(0)} V_{(\beta)} (\divg k - d(\tr k)) \\
& \quad = 2 \iota_{V_{(\beta)} \nabla V_{(0)} - V_{(0)} \nabla V_{(\beta)}} \left(k - (\tr k) b\right)+ 2 V_{(0)} V_{(\beta)} \left( \divg_g k - d (\tr_g k) + \mathcal{O}(e^{-(2\tau + 1)r}) \right) \\
& \quad = 2 \iota_{V_{(\beta)} \nabla V_{(0)} - V_{(0)} \nabla V_{(\beta)}} \left(k - (\tr k) b\right)+ 2 V_{(0)} V_{(\beta)} J 
+ \mathcal{O}(e^{-(2\tau - 1)r})  .
\end{split}
\end{align*}
Here and in the rest of the proof all quantities are computed with respect to hyperbolic metric $b$, unless otherwise indicated. Using the Poincar\'e ball model for the hyperbolic metric as in Section \ref{Charges}, it is straightforward to check that
\begin{align*}\labeleq{eqCi}
\begin{split}
 V_{(i)} \nabla V_{(0)} - V_{(0)} \nabla V_{(i)} & =  V^2_{i} \,\nabla \left(\frac{V_{(0)}}{V_{(i)}}\right) \\
 & = \frac{4 (y^i)^2}{(1-\vert y\vert ^2)^2} \,\nabla \left(\frac{1+\vert y\vert ^2}{2 y^i}\right)\\
 & = (y^i)^2 \sum_{j=1}^n \left(\frac{1+\vert y\vert ^2}{2 y^i}\right)'_{y^j} \partial_{y^j}\\
 & =  y^i y^j \partial_{y^j} - \frac{1+\vert y\vert ^2}{2} \partial_{y^i}\\
 & =  C_{(i)}
\end{split}
\end{align*}
for all $i = 1, \ldots, n$. Hence 
\begin{align*}
 \frac{d}{dt}{\vert _{t=0}}\, \bU_{(V_{(i)},0)} (g,k) &= \bU_{(0,C_{(i)})} (g,k)+ 2 V_{(0)} V_{(i)} J + \mathcal{O}(e^{-(2\tau - 1)r})  
\end{align*}
and
\begin{align*}
 \frac{d}{dt}{\vert _{t=0}} \,\bU_{(V_{(0)},0)} (g,k) &= 2 V_{(0)}^2 J + \mathcal{O}(e^{-(2\tau - 1)r}).  
\end{align*}
The claim follows from the asymptotics of $J$.
\qed

\begin{remark}
In fact, Theorem \ref{thEvolution} applies under more general assumptions on lapse and shift, in particular when $N=V_{(0)}( 1 + \mathcal{O}_1(e^{-\frac{n}{2} r}))$ and $X=\mathcal{O}_2(e^{-(n+\varepsilon) r})$ for some $\varepsilon>0$ as $r\to\infty$. The details are left to the reader.
\end{remark}

\begin{remark}
In view of Theorem \ref{thEvolution} and the Newtonian/special relativity law $\frac{d}{dt} (m \vec{z}) = \vec{P}$ also alluded to in Section \ref{sect:AE}, we propose to call $\mathbb{Q}_{(0,C_{(i)})}$ the \emph{linear momentum} of the initial data set $(M,g,k)$, thereby again deviating from the definition by Chru\'sciel, Maerten, and Tod \cite{ChruscielMaertenTod}, see \eqref{eqMassVector} and \eqref{eqCoM}. 
\end{remark}

\section{Conclusion}
We have defined a new notion of center of mass for asymptotically hyperbolic initial data sets as in Definition \ref{defAdSData}, see Definition \ref{defCoMcharges}. This center of mass is a point in the hyperbolic reference space, coinciding with the limit of the hyperbolic coordinate centers of the leaves of the CMC-foliation constructed by Neves and Tian and extended to non-balanced coordinates in Theorem \ref{thCoMUnbalanced}, see Section \ref{sect:geometric}. The new center transforms adequately under hyperbolic isometric changes of the asymptotic coordinates by definition and evolves correctly under the Einstein evolution equations with appropriate cosmological constant (Theorem \ref{thEvolution}).

\vfill
\pagebreak

\bibliographystyle{amsplain}
\bibliography{biblio-AH-mass}

\providecommand{\bysame}{\leavevmode\hbox to3em{\hrulefill}\thinspace}
\providecommand{\MR}{\relax\ifhmode\unskip\space\fi MR }
% \MRhref is called by the amsart/book/proc definition of \MR.
\providecommand{\MRhref}[2]{%
  \href{http://www.ams.org/mathscinet-getitem?mr=#1}{#2}
}
\providecommand{\href}[2]{#2}
\begin{thebibliography}{10}

\bibitem{AnderssonCaiGalloway}
Lars Andersson, Mingliang Cai, and Gregory~J. Galloway, \emph{Rigidity and
  positivity of mass for asymptotically hyperbolic manifolds}, Ann. Henri
  Poincar\'e \textbf{9} (2008), no.~1, 1--33.

\bibitem{ADM}
Richard Arnowitt, Stanley Deser, and Charles~W. Misner, \emph{{Coordinate
  Invariance and Energy Expressions in General Relativity}}, Phys. Rev.
  \textbf{122} (1961), no.~3, 997--1006.

\bibitem{Bartnik}
Robert Bartnik, \emph{{The Mass of an Asymptotically Flat Manifold}},
  Communications on Pure and Applied Mathematics \textbf{39} (1986), 661--693.

\bibitem{BartnikIsenberg}
Robert {Bartnik} and James {Isenberg}, \emph{{The constraint equations.}}, {The
  Einstein equations and the large scale behavior of gravitational fields. 50
  Years of the Cauchy problem in general relativity. With DVD}, Basel:
  Birkh\"auser, 2004, pp.~1--38 (English).

\bibitem{BO}
Robert Beig and Niall~{\`O} Murchadha, \emph{The poincar\'e group as the
  symmetry group of canonical general relativity}, Ann. Physics \textbf{174}
  (1987), no.~2, 463--498.

\bibitem{Cederbaum_newtonian_limit}
Carla Cederbaum, \emph{{The Newtonian Limit of Geometrostatics}}, Ph.D. thesis,
  FU Berlin, 2012,
  \href{http://arxiv.org/abs/1201.5433v1}{arXiv:\linebreak[1]1201.5433v1}.

\bibitem{Cederbaum-Nerz}
Carla Cederbaum and Christopher Nerz, \emph{{Explicit Riemannian manifolds with
  unexpectedly behaving center of mass}}, Ann. Henri Poincar\'e (2014), online
  first: http://tiny.cc/cttcox.

\bibitem{CWYEuclidean}
Po-Ning Chen, Mu-Tao Wang, and Shing-Tung Yau, \emph{Conserved quantities in
  general relativity: from the quasi-local level to spatial infinity},
  arxiv:1312.0985v1.

\bibitem{ChenWangYauQL}
Po-Ning Chen, Mu-Tao Wang, and Shing-Tung Yau, \emph{Quasilocal angular
  momentum and center of mass in general relativity},
  \href{http://arxiv.org/abs/1312.0990v1}{arXiv:\linebreak[1]1312.0990v1},
  2013.

\bibitem{ChenWangYauhyperbolic}
\bysame, personal communication, 2014.

\bibitem{ChristodoulouYau}
Demetrios Christodoulou and Shing-Tung Yau, \emph{Some remarks on the
  quasi-local mass}, Mathematics and general relativity ({S}anta {C}ruz, {CA},
  1986), Contemp. Math., vol.~71, Amer. Math. Soc., Providence, RI, 1988,
  pp.~9--14.

\bibitem{Chrus2}
Piotr~T. Chru{\'s}ciel, \emph{{On the invariant mass conjecture in general
  relativity}}, Communications in Mathematical Physics \textbf{120} (1988),
  233--248.

\bibitem{ChruscielHerzlich}
Piotr~T. Chru{\'s}ciel and Marc Herzlich, \emph{The mass of asymptotically
  hyperbolic {R}iemannian manifolds}, Pacific J. Math. \textbf{212} (2003),
  no.~2, 231--264.

\bibitem{ChruscielMaertenTod}
Piotr~T. Chru{\'s}ciel, Daniel Maerten, and Paul Tod, \emph{{Rigid upper bounds
  for the angular momentum and centre of mass of non-singular asymptotically
  anti-de Sitter space-times}}, Journal of High Energy Physics \textbf{2006}
  (2006), no.~11, 084.

\bibitem{ChruscielNagy}
Piotr~T. Chru{\'s}ciel and Gabriel Nagy, \emph{The mass of spacelike
  hypersurfaces in asymptotically anti-de {S}itter space-times}, Adv. Theor.
  Math. Phys. \textbf{5} (2001), no.~4, 697--754.

\bibitem{Cortier}
Julien {Cortier}, \emph{{A family of asymptotically hyperbolic manifolds with
  arbitrary energy-momentum vectors.}}, {J. Math. Phys.} \textbf{53} (2012),
  no.~10, 102504, 15.

\bibitem{CDG}
Julien Cortier, Mattias Dahl, and Romain Gicquaud, \emph{Mass-like invariants
  for asymptotically hyperbolic metrics}, in preparation.

\bibitem{EichmairMetzgerInventiones}
Michael Eichmair and Jan Metzger, \emph{{Unique isoperimetric foliations of
  asymptotically flat manifolds in all dimensions}}, Invent. Math. (2013),
  arXiv:1204.6065v1, to appear.

\bibitem{Lan_Hsuan_Huang__Foliations_by_Stable_Spheres_with_Constant_Mean_Curvature}
Lan-Hsuan Huang, \emph{Foliations by {S}table {S}pheres with {C}onstant {M}ean
  {C}urvature for isolated systems with general asymptotics}, Communications in
  Mathematical Physics \textbf{300} (2010), no.~2, 331--373.

\bibitem{HuiskenYau}
Gerhard Huisken and Shing-Tung Yau, \emph{{Definition of center of mass for
  isolated physical systems and unique foliations by stable spheres with
  constant mean curvature.}}, {Invent. Math.} \textbf{124} (1996), no.~1-3,
  281--311.

\bibitem{Maerten}
Daniel Maerten, \emph{Positive energy-momentum theorem for
  {A}d{S}-asymptotically hyperbolic manifolds}, Ann. Henri Poincar\'e
  \textbf{7} (2006), no.~5, 975--1011.

\bibitem{MazzeoPacard}
Rafe Mazzeo and Frank Pacard, \emph{{Constant curvature foliations in
  asymptotically hyperbolic spaces.}}, {Rev. Mat. Iberoam.} \textbf{27} (2011),
  no.~1, 303--333 (English).

\bibitem{MichelMass}
Beno\^it Michel, \emph{Geometric invariance of mass-like asymptotic
  invariants}, J. Math. Phys. \textbf{52} (2011), no.~5, 052504, 14.

\bibitem{Moller}
Christian M{\o}ller, \emph{On the definition of the centre of gravity in an
  arbitrary closed system in the theory of relativity}, A, Theoretical physics,
  vol.~3, Commun. Dublin Inst. Advanced Studies, 1949.

\bibitem{nerz2013timeevolutionofCMC}
Christopher Nerz, \emph{Time evolution of {ADM} and {CMC} center of mass in
  general relativity},
  \href{http://arxiv.org/abs/1312.6274v1}{arXiv:\linebreak[1]1312.6274v1},
  2013.

\bibitem{Nerz-characterization}
\bysame, \emph{{Geometric characterizations of asymptotic flatness and linear
  momentum in general relativity}},
  \href{http://arxiv.org/abs/1409.6039v1}{arXiv:\linebreak[1]1409.6039v1},
  2014.

\bibitem{NevesTian1}
Andr\'e Neves and Gang Tian, \emph{{Existence and uniqueness of constant mean
  curvature foliation of asymptotically hyperbolic 3-manifolds.}}, {Geom.
  Funct. Anal.} \textbf{19} (2009), no.~3, 910--942.

\bibitem{NevesTian2}
\bysame, \emph{{Existence and uniqueness of constant mean curvature foliation
  of asymptotically hyperbolic 3-manifolds. II.}}, {J. Reine Angew. Math.}
  \textbf{641} (2010), 69--93.

\bibitem{RT}
Tullio Regge and Claudio Teitelboim, \emph{Role of surface integrals in the
  hamiltonian formulation of general relativity}, Annals of Physics \textbf{88}
  (1974), no.~1, 286--318.

\bibitem{Rigger}
Ralf Rigger, \emph{{The foliation of asymptotically hyperbolic manifolds by
  surfaces of constant mean curvature (including the evolution equations and
  estimates).}}, {Manuscr. Math.} \textbf{113} (2004), no.~4, 403--421.

\bibitem{Szabados-old}
L\'aszl\'o Szabados, \emph{On the {P}oincar\'e {s}tructure of {a}symptotically
  {f}lat {s}pacetimes}, Class. Quantum Grav. \textbf{20} (2003), 2627--2661.

\bibitem{Szabados}
\bysame, \emph{The {P}oincar\'e {S}tructure and the {C}entre-of-{M}ass of
  {A}symptotically {F}lat {S}pacetimes}, Analytical and Numerical Approaches to
  Mathematical Relativity, Lect. Notes Phys., vol. 692, Springer, 2006,
  pp.~157--184.

\bibitem{WangMass}
Xiaodong Wang, \emph{The mass of asymptotically hyperbolic manifolds}, J.
  Differential Geom. \textbf{57} (2001), no.~2, 273--299.

\end{thebibliography}
\vfill
\end{document}